\documentclass{amsart}[14pt]

\usepackage{amsmath, amstext, amssymb, amsopn, amsthm, amscd, amsfonts}
\usepackage{graphicx}
\usepackage{epsfig}

\newcommand{\R}{\mathbb R}

\newcommand{\C}{\mathbb C}

\newcommand{\HP}{\mathbb H}

\newtheorem{theo}{Theorem}

\newtheorem{lema}{Lemma}

\newtheorem{prop}{Proposition}

\newtheorem{defi}{Definition}

\newtheorem{rema}{Remark}

\begin{document}

\title{Univalent Baker domains and boundary of deformations}         
\author{P. Dom{\'i}nguez and G. Sienra}      
\date{}         
\maketitle

\begin{abstract} For $f$ an entire transcendental map with a univalent Baker domain $U$ of hyperbolic type, we study pinching deformations in $U$, the support of this deformation being certain laminations $\mathcal {R}(\Lambda)$ in the grand orbit of $U$. We show that if $f$ satisfies certain conditions, the limit of the deformation is an entire transcendental map $F$. Such map contains a family  of simply connected wandering domains, disjoint from the postcritical set. We interpret these results in terms of the Teichm\"uller space of $f$, $Teich(f)$ as subset of ${\mathcal{M}}_{f}$, the space of marked topologically equivalent maps to $f$.



\end{abstract}

\section{introduction}

This paper studies some properties of the dynamics of entire transcendental maps which relates univalent Baker domains and wandering domains. 

The first example of an entire transcendental function with such domains was given by Fatou [F]. Examples of Baker domains of higher period were given by Baker, Kotus and L\"u [BKL].  For complete references see the survey of Bergweiler [B1]. Existence of Baker domains without singularities in its interior (univalent Baker domains) were constructed by Bergweiler in [B2], Rippon and Stallard [RS], Baransky and Fagella [BF], Fagella and Hendriksen [FH1] and other mathematicians. In section 2 will be given the notation, definitions and concepts that will be useful through the paper.

In this article we are interested in the limit of deformations of entire transcendental maps $f$ with a periodic univalent Baker domain $U$. Specifically our domain $U$ is a Baker domain of hyperbolic type I  according to [BF]. We apply pinching deformations along curves in these domains.  Pinching deformation was introduced by Makienko [M] in the context of rational maps in analogy of Kleinian groups, subsequently [T], [HT], [BT] and others, developed and applied the theory.

A periodic univalent Baker domain $U$ admits an hyperbolic metric, thus geodesic curves are well defined. In Section 3 we define and study certain geodesic laminations $\Lambda$ in $U$, that we call Baker laminations. The grand orbit of such laminations, $\mathcal {R}$ $(\Lambda)$, is where the pinch deformation occurs. These laminations are invariant under the action of the map $f$. In section 4, we define pinching deformations.

In section 5 we study a special kind of laminations. Those that contains the union of geodesics that join the expansive periodic points in the boundary of $U$ with $\infty$, here denoted by $\lambda_{\infty}$.  We prove in Theorem 1, that for any lamination that contains the geodesics $\lambda_{\infty}$, pinching deformation is divergent.

In the other hand in section 6 we prove:

{\bf Theorem 2.} {\sl Let $f$ be an entire transcendental map satisfying:

(a) $f$ has a univalent Baker domain of hyperbolic type I,

(b) the postcritical set, is a positive distance away from the Julia set and

(c) $J(f)$ is thin at $\infty$.

Then for  $\Lambda$ any Baker lamination in $U$ which does not contains $\lambda_{\infty}$, the pinching process along $\mathcal {R}(\Lambda)$ converges uniformly to an entire transcendental map $F$ which exhibits in its Fatou set a family of bounded simply connected wandering domains, disjoint from the postcritical set}.

The hypothesis of the theorem allow the map $f$ to be expanding in the Julia set according to a theorem of Stallard ([St], Theorem A) and so it is weakly hyperbolic (Section 2.3), hence the convergence techniques in the paper of Haissisky and Tan Lei [HT] can be used in this case. This is the main reason to have chosen such kind of Baker domains.

In section 7, we study pinching deformations as seeing in the marked Hurwitz space of the map $f$, here denoted by $M_{f}$ and which is defined as the space of all marked maps topologically equivalent to $f$ (Definition 4). We include the Teichm\"uller space of $f$, $Teich(f)$ into $M_{f}$. Theorem 1 shows that pinching along $\lambda_{\infty}$ does not converges to a point in $M_{f}$, hence $Teich(f)$ is not compactly contained in $M_{f}$. Pinching along the other Baker laminations defines (by Theorem 2) points at the boundary of $Teich(f)$ contained in $M_{f}$. 


Finally, it is shown that if we take a sequence of deformations in the slice $Teich(U/f)$ of $Teich(f)$, with the property that it converges uniformly to an entire transcendental map at the boundary of $Teich(f)$ (as subset of $M_{f}$) then, the only development of $U$ towards a periodic Fatou domain (in case it exists) is to a new periodic univalent Baker domain.

\section{Univalent Baker domains}

\subsection{Notation and definitions}\

Our notation shall be: $\HP$ for the upper half plane, ${\R}^{*}={\R} \cup \{\infty\}$ is the boundary of the upper half plane, $\C$ the complex plane, $\hat{\C}$ the Riemann sphere, $d$ the euclidian distance in $\C$, $d_{h}(W)$ the hyperbolic metric in a simply connected domain $W \subset {\C}$ and $l_{W}(s)$ the hyperbolic length of an arch $s$ in the hyperbolic metric of $W$. If $U$ is a set, we write $U/f$ for the quotient of $U$ by the action of $f$, $f'(w)$ is the derivative of $f$ at $w$.

The map $f$ defines a partition of the plane into two sets, the Fatou set $F(f)$ and the Julia set $J(f)$. The Fatou set is where the iterates of $f$ form a normal family and the Julia set is the complement of the Fatou set. Properties of these sets can be found in [B1].

The singularities of the inverse function are the critical values and the asymptotic values. Critical values are the images of points $c$ such that $f'(c)=0$, asymptotic values are points $a \in \C$ for which there exist a path ${\gamma}(t) \rightarrow \infty$ as $t \rightarrow \infty$ such that $f({\gamma}(t))\rightarrow a$ as $t \rightarrow \infty$. We denote by $Sing(f^{-1})$ the closure of the set of critical values and asymptotic values and by $P(f)$, the postcritical set, which is the closure of the (positive) orbits of points in $Sing(f^{-1})$.

\subsection{Univalent Baker domains}\

The classification of Fatou domains for transcendental maps are well known, see for instance [B1]. Fatou domains that do not appear for Rational maps are the following.

\begin{defi}({\bf Baker domain}) For $f$ an entire transcendental map, $U$ a periodic component of the Fatou set of period $p$, is called a Baker domain if for all $z \in U$, $f^{np}(z) \rightarrow {\infty}$ as $n \rightarrow \infty$.
\end{defi}

It was proved by Baker in [Ba] that all Baker domains for entire transcendental maps are simply connected. In [EL] it is established that if $Sing(f^{-1})$ is bounded then $f$ has not Baker domains. It is not necessarily for a Baker domain to have any of the singular values inside. 

A {\bf Univalent Baker domain} is a periodic Baker domains on which $f^{p}$ is univalent and these are the kind of domains that we are interested in this paper. For entire transcendental maps, a classification of this domains is given by Baransky, Hendriksen and Fagella [BF], [FH1] as follows. 

Assume that $f$ has an invariant Baker domain $U$, then there exist a point $\xi \in \hat{\C}$ such that the backward iterates under $(f|_{U})^{-1}$ of all points in $U$ tend to $\xi$ through the same access (they call it the backward dynamical access), and one of the following cases occurs:

(a) If $\xi \neq \infty$ is a fixed point in the boundary of $U$, then $U$ is of {\sl hyperbolic type I}.

(b) If $\xi=\infty$ but the backward dynamical access is different than the forward dynamical access, then $U$ is of {\sl hyperbolic type II}.

(c) If $\xi=\infty$ but the backward dynamical access is the same to the forward one, then $U$ is of {\sl parabolic type}.

The above domains have a uniformization $\Psi:\HP \rightarrow U$ and a linear map $G:{\HP} \rightarrow {\HP}$ such that $(f|_{U}) \circ \Psi=\Psi \circ G$. If $U$ is of hyperbolic type I or II, then $G(z)=az$, $a > 1$. For $U$ of parabolic type, then $G(z)=z+1$. Thus the diagram bellow commutes.

\begin{center}
\hspace{0cm} $f$ 

$U \hspace{.5cm} \longrightarrow \hspace{.5cm} U$ 
\end{center}

\begin{equation}
\hspace{2cm} {\Psi} \Big\uparrow \hspace{2cm} \Big\uparrow {\Psi} \hspace{2cm} 
\end{equation}

\begin{center}
\hspace{0cm} $G$

${\HP} \hspace{.5cm} \longrightarrow \hspace{.5cm} {\HP}$

\end{center}

Now, let us denote by $\alpha$ the geodesic $it \in {\HP}$, $t>0$ and by $A_{a}$ the annulus $\HP/G$ with modulus $m(A_{a})={\pi}/log(a)$. Let $c={\alpha}/G$ be the core geodesic of $A_{a}$, with length $log(a)$. The dynamic action of $f$ in $U$ defines the quotient space $U/f$.  Hyperbolic domains of type I and II, have as quotient space $U/f$, an annulus $A_{a}$ which is conformally the same that the quotient annulus $\HP/G$.

\subsection{Hyperbolicity}\

In this paper we will be interested in entire transcendental functions with a univalent periodic domain such that has the postcritical set at a positive distance away from the Julia set or equivalently $d(P(f),J(f))>C>0$. This kind of functions are non-uniformly hyperbolic according to Theorem A in [St], in the sense that for all $z\in J(f)$, $|f^{n'}(z)| \rightarrow \infty$ as $n$ tends to $\infty$. However the condition that is relevant for us is the following.

\begin{defi}
({\bf Weak hyperbolicity}) We say that $f$ is weakly hyperbolic if there are constants $r>0$ and ${\delta}< \infty$ such that, for all $z \in J_{f}-\{preparabolic$ $points\}$, there is a subsequence of iterates $\{f^{n_k}\}$ such that

\begin{displaymath}
deg(f^{n_k}:W_{k}(z) \rightarrow D(f^{n_k}(z),r)) \leq {\delta},
\end{displaymath}

where $W_{k}(z)$ is the connected component of $f^{-n_{k}}(D(f^{n_k}(z),r))$ containing $z$.
\end{defi}

Due to the condition of the postcritical set and the Julia set above, our maps do not have parabolic points.

\subsection{Examples of hyperbolic functions with univalent Baker domains}\

 We are interested in entire transcendental weakly hyperbolic maps with $d(P(f),J(f))>C>0$ and with a univalent Baker domain, then the Bergweiler's example is the typical one.

Consider the map $f(z)=2-log2+2z-e^{z}$. It was shown in [B2] that $f$ has an invariant Baker domain $U$ containing the half plane $\{Re(z)<-2\}$ and the boundary of $U$ is a Jordan curve, all singular points are in Fatou domains and the plane distance from $\partial U$ to the postcritical set is positive. The Singular points of $f$ are located in one superattracting component and in a union of simply connected wandering domains. Moreover,
 the Lebesgue measure of $J(f)$ is zero, as it is shown in [St].

Following the example of Bergweiler, we can define, for $m>1$, the family of maps.

\begin{displaymath}
f_{m}(z)=(m-1)log(\frac{2}{2m-1})+\frac{2m-1}{2}+mz-e^{z}
\end{displaymath}

The family $f_{m}(z)$ has the same properties that $f$ except that we have changed the superattracting fixed point by an attracting fixed point at $\frac{2m-1}{2}$. For each of these functions, the map $G$ of diagram (1) is $G(z)=mz$. 

Since the maps $f_{m}(z)$ satisfy the condition $d({P(f)},J(f))>C>0$, they are hyperbolic (and weakly hyperbolic).

\section{Baker laminations}

Pinching deformation is supported in geodesic arcs with good neighborhoods, defined as follows.

We consider a band $B_{\delta}=[{\pi}/2-{\delta},{\pi}/2+{\delta}]{\times}{\R}$ in the complex plane, and consider
its image $V_{\delta}$ under the exponential map, with $\delta$
so $V_{\delta} \subset {\HP}$. Thus $V_{\delta}$ is a {\sl good neighborhood of thickness ${\delta}$} of the horizontal geodesic ${\alpha}$. If ${\gamma}$ is any other complete
 geodesic in ${\HP}$, then there is an isometry $A$ of the upper half plane such that $A({\alpha})={\gamma}$; we say that $A(V_{\delta}):=V_{\delta}(\gamma)$ is a good neighborhood (of thickness $\delta$) for ${\gamma}$. Here the line $({\pi}/2){\times}{\R}$ is sent to the geodesic $\gamma$ under $A \circ Exp(z)$.

From now on, $U$ {\bf denotes a univalent Baker domain}. Consider a geodesic $\lambda \in U$, we have that ${\Psi}^{-1}{\lambda}=\gamma$ is a geodesic in $\HP$ with a good neighborhood $V_{\delta}(\gamma)$. We define                            ${\Psi}V_{\delta}(\gamma):=V_{\delta}(\lambda)$ as a {\bf good neighborhood} of $\lambda$ of thickness $\delta.$

\begin{defi}

Let $\Lambda$ be a set of (complete and open) geodesics in $U$, we say that $\Lambda$ is a {\bf Baker lamination in $U$}, if the elements in $\Lambda$ satisfies:

(i) It is invariant under the action of $f$. That is, if $\lambda \in \Lambda$, then $f^{n}(\lambda) \in \Lambda$ for any $n$ a positive integer. Also if we consider the inverse branch of $f^{-n}$ from $U$ to $U$, we have that $f^{-n}(\lambda) \in \Lambda$.

(ii) For any $\lambda,\lambda'\in \Lambda$, then $\lambda \bigcap \lambda'=\emptyset.$


(iii) Its end points are well defined points in $\partial U$ and any two geodesics have different end points in $\partial U$.


(iv) There are not self accumulating leaves. That is, all geodesics in $\Lambda$ admit a good neighborhood of thickness $\delta$ with the property that good neighborhoods are pairwise disjoint.

\end{defi}

See figure 2 (page 11), for an example.

Let us denote by $End(\gamma)$ the two end points of the geodesic $\gamma$ and by $\bar{\gamma}={\gamma} \cup End(\gamma)$. We observe that condition (iii) can be restated as $End(\gamma) \bigcap End(\gamma')=\emptyset$ for any ${\gamma},{\gamma}' \in \Lambda$.

If the Baker domain is periodic with $U={\bigcup}_{i=1}^{n-1}U_{i}$ and $fU_{i}=U_{i+1}$, then 
$f^{n}U_{1}=U_{1}$ and we can consider a Baker lamination $\Lambda_{1}$ in $U_{1}$, invariant under $f^{n}$ and so $f^{k}{\Lambda_{1}}={\Lambda}_{k}$, $k=1,...,n-1$ define a Baker lamination in $U_{k+1}$ respectively, invariant under $f^{n}$. The construction  defines to  $\Lambda=\bigcup \Lambda_{i=1}^{n}$ as a Baker lamination in the periodic domain $U$. Good neighborhood in $\Lambda$ are constructed similarly, first for the Baker lamination in $U_{1}$ and then by iteration for any element of $\Lambda$.

{\bf Notation}: We will denote by $\mathcal {R}(\Lambda)={\cup}_{k}f^{-k}(\Lambda)$, the full orbit of $\Lambda$. Let ${\overline{\Lambda}}$ be the set of geodesics in $\Lambda$ together with its end points and denote by ${\mathcal {R}}(\overline{\Lambda})={\cup}_{k}f^{-k}(\overline{\Lambda})$, the full orbit of ${\overline{\Lambda}}$  . Let us denote by $\mathcal{V}$ $({\Lambda})$ $={{\bigcup}_{k}}f^{-k}(V_{\delta}({\Psi}{\Lambda}))$, the full orbit of the good neighborhoods. Let $\mathcal{U}$ the full orbit of the Baker domain $U$.

\subsection{Paths in $\mathcal {R}(\Lambda)$}\

 The main aim in this section is to show that there are not two geodesics of a Baker lamination that touch at a common point.

We know from [BE] that the inverse images of $U$ are disconnected. Let $U_{1}$ one of such components, then if there is a singularity of $f$ in $U_{1}$ then $f:U_{1} \rightarrow U$ has degree $d \geq 2$, hence there is a point $w$ in $\partial U_{1}$ such that $f(w)=\infty$. Therefore, no preimage contains critical points and $f^{-1}(U)$ is made of infinitely many components. Recursively, this implies that no critical point is in the full orbit of $U$, hence for any $V \subset  {\mathcal{U}}$ we have that
$V$ is connected , in fact it is simply connected (Theorem 9 in [B1]) and $f^{n}:V \rightarrow U$ is one to one. 

The problem now, is to know if there are periodic Baker domains or inverse images that may have some boundary points in common (other than $\infty$). In [Si] there is a construction of functions with invariant univalent Baker domains $U$ and $U'$ (or more), which touch each other at their common fixed point $p=\partial U \cap \partial U'$. The following lemma implies that this is the only possible case.

Assume that $f$ has a family of periodic univalent Baker domains $\{U_{i}\}$ and $d(Sing(f^{-1}),J(f))>0$. Let $\Lambda$ a Baker lamination in $U_{i}$, then we have:


\begin{lema}
If $\Lambda$ does not contains $\lambda_{\infty}$ and $\gamma_{1}$, $\gamma_{2}$ are in $\mathcal {R}(\Lambda)$, then  $\bar{\gamma_{1}} \bigcap \bar{\gamma_{2}}= \emptyset$
\end{lema}

\begin{proof}Suppose that there are $\gamma_{1}$ and $\gamma_{2}$ such that $\bar{\gamma_{1}} \bigcap \bar{\gamma_{2}}=x$. Let us consider two cases.

Case (a). Let $U$, $V$ be fixed univalent Baker domains and assume that $\partial U \bigcap \partial V = x$. Then $f^{n}(x)=x_{n}$ is also in $\partial U \bigcap \partial V$. Consider paths $\sigma_{1} \subset U$ and $\sigma_{2} \subset V$ from $x$ to $x_{1}$ respectively and denote by $D$ the disc which is the interior of  $\sigma=\sigma_{1} \cup \sigma_{2}$. Then $f(D)$ takes its maximum over $\sigma$, therefore the disc $D_{1}$ interior to $f(\gamma)$ satisfies $f(D)=D_{1}$. Consequently, when $n$ tends to infinity, $f^{n}(\sigma)$ tends to $\infty$, so $f^{n}(D)$ tends to $\infty$, which is a contradiction since $D$ contains points of the Julia set. 

Case (b). Let $U$ be a univalent Baker domain such that $f(U)=U$ and $V \subset {\mathcal {U}}$, then there exist $n>0$ with $f^{n}(V)=U$, if $\partial U \bigcap \partial V = x$. Let $g=f^{n}$, then $g$ sends a neighborhood of $x$ into its image in a two to one fashion, then $x$ is a singular point in $J(g)=J(f^{n})=J(f)$ a contradiction with the assumption that $d(Sing(f^{-1}),J(f))>0$.

The general situation is proved by iteration of $\gamma_{1}$ and $\gamma_{2}$ until we fall in Case (a) or (b).

\end{proof}

\begin{rema} From Lemma 1, we have that if $\tau$ is a connected path in ${\mathcal {R}(\overline{\Lambda})}$, then between any to geodesics in $\tau$ there is always another geodesic. This path is made up by a set of geodesics belonging to the great orbit of the lamination $\Lambda$ and a Cantor set which is the set of boundary points of laminations or the limit points of sequences of geodesics in $\tau$. This Cantor set belong to the Julia set. In the process of pinching, these geodesics tend to be shorter and shorter until their two ends are identified to a point. At the end of this process the path $\tau$ is shrinked into a path but not to a point.
\end{rema}

\subsection{The $\lambda_{\infty}$-geodesic}\

Consider the geodesic $\alpha$ in $\HP$ of Section 2.2 and $\Psi$ the uniformization of some univalent Baker domain $U$ of hyperbolic type I. Then $\Psi({\alpha})$ is a lamination in $U$, that we will denote by $\lambda_{\infty}$, see figure 1. Its projection to $U/f$ is the core geodesic $c$. If $f|U$ is of hyperbolic type I, the geodesic $\lambda_{\infty}$  joins the fixed point $p \in \partial U$ with $\infty$ and it is invariant under the action of $f$, hence a Baker lamination. In the case of a periodic univalent Baker domain $U_{1},...,U_{k}$, there is a periodic expanding point $q_{1},...,q_{k}$ with $q_{i} \in \partial U_{i}$ and therefore we will consider $\lambda_{\infty}$ to be the union of the geodesics that joins $q_{i}$ with $\infty$.

\begin{figure}[h]
\begin{center}
\includegraphics[scale=0.5]{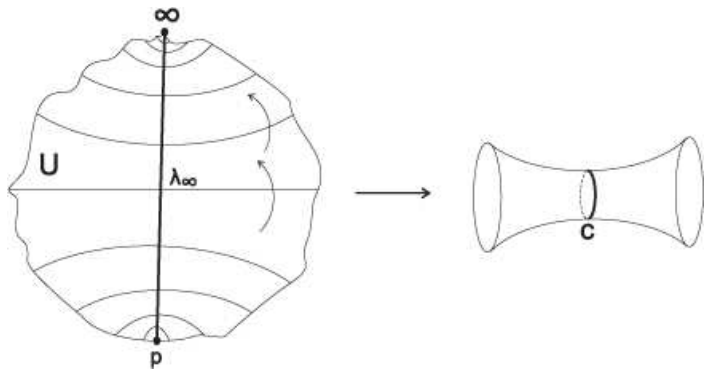}
\caption{A Baker domain $U$ with the geodesic ${\lambda}_{\infty}$ and the quotient by the action $U/f$, an annulus. The closed geodesic $c$ is the image of ${\lambda}_{\infty}$.}
\end{center}
\end{figure}

Observe that a geodesic in $U$ that contains $\infty$ as end point, projects in $U/f$ to a geodesic such that necessarily accumulates on the core geodesic $c$. Consequently, by property (iv) we have that if a leaf of a Baker lamination has $\infty$ as an end point, it has to be $\lambda_{\infty}$.

\subsection{General Baker laminations}\

Now consider $U$ a forward invariant univalent Baker domain and two points $u,v$ in ${\partial U}-\{{\infty}\}$ which are accessible from $U$. These points always exist and dynamically satisfy that $f^{n}(u) \rightarrow \infty$, $f^{n}(v) \rightarrow \infty$ as $n \rightarrow \infty$. Consider a geodesic $\gamma_{u,v}$ in $U$ joining $u$ and $v$. Choose the points $u,v$ in such a way that $\gamma_{u,v} \bigcap f(\gamma_{u,v})=\emptyset$, then the forward an backward orbit of $\gamma_{u,v}$ in $U$ is a Baker lamination.

In this way it is possible to consider finite couples of points $u_{i}, v_{i}$ in ${\partial U}-\{{\infty}\}$ and accessible from $U$. If the geodesics joining them are disjoint as well as their images, then their forward an backward orbits in $U$ under $f$ are a Baker lamination. 

We can add to such laminations the geodesic $\lambda_{\infty}$ if the set satisfies (i) to (v), then we obtain laminations that are consider in Theorem 1.

\section{pinching deformation}

We assume that $f$ is an entire transcendental map satisfying the following conditions:

(a) It has at least one periodic univalent Baker domain $U_{1},...,U_{n}$ of hyperbolic type I.

(b) The postcritical set is a positive distance away from the Julia set.

Under the above conditions $f$ is weakly hyperbolic where we can choose ${\delta}=1$ and $r=(1/2)dis(J(f)$, Post($f$)). The examples of in Section 2.1. satisfies these conditions.

 A {\bf pinching combinatorics} is defined as follows:

${\bullet}$ for $\delta >0$, a set of finite Baker laminations of thickness $\delta$, ${\Lambda}_{1} \subset U_{1},...,{\Lambda}_{n} \subset U_{n}$ such that $f^{k}{\Lambda}_{i}={\Lambda}_{j}$, $k \geq 0$, $i,j=1,...,n$. Observe that each ${\lambda_{i}} \in {\Lambda}$ is in $F(f)$ and disjoint from the orbits of the critical points.

${\bullet}$ the $\lambda_{i}$'s are mutually disjoint.

${\bullet}$ $f:{\lambda_{i}} \rightarrow {\lambda_{j}}$ is a homeomorphism.

${\bullet}$ An invariant set $\mathcal{V}$ $({\Lambda})$, of disjoint good neighborhoods for the Baker lamination and in consequence for $\mathcal{R}({\Lambda})$.

For a {\bf pinching deformation}, we take from [HT] the quasiconformal symmetric model on strips, with our $B_{\delta}$ as the support of the model strip. In that paper is defined the quasiconformal map ${\tilde{P}}_{t}(x+iy)=x+iv_{t}(y)$, which has the following four properties:

(1) It commutes with the translation by any real number.

(2) It is the identity in the sub-strip $\{-L_{y} \leq y \leq L_{y}\}$, being $L_{y}$ as in [HP].

(3) The coefficient of the Beltrami form

\begin{displaymath}
\frac{\partial\tilde{P}/\partial \bar{z}}{\partial\tilde{P}/\partial z}\vert_{x+iy}=\frac{1-
\frac{\partial}{\partial y} v_{t}(y)}{1+\frac{\partial}{\partial y} v_{t}(y)}
\end{displaymath}

is continuous in $(t,x+iy)\in[0,1] \times$ $\{-L_{r} \leq y \leq L_{r}\}$, whose norm is locally uniformly bounded from 1 if $(t,x) \neq (1,L_{r})$  and tends to 1 as $(t,y) \rightarrow (1,L_{r})$.

\medskip

If $V$ is a connected component of $\mathcal{V}({\Lambda})$, define the conformal map  $\psi:V \rightarrow B_{\delta}$ as the well defined inverse branch of the map ${\Psi} \circ A \circ Exp:B_{\delta} \rightarrow U_{i} $ that sends the strip to a good neighborhood of an element of the grand orbit of a Baker lamination. For $t \in [0,1[$, set $\sigma'_{t}=(\tilde{P} \circ \psi)^{*}(\sigma_{0})$ to the pull back of the standard complex structure on $B_{\delta}$. We spread $\sigma'_{t}$ to the whole orbit in $\mathcal{V}({\Lambda})$ under the map $f$. Let $\sigma_{t}$ be the extension to the Riemann sphere by setting $\sigma_{t}=\sigma_{0}$ in the complement of $\mathcal{V}({\Lambda})$. Then, $\sigma_{t}$ is an f-invariant complex structure.

The family $({\sigma}_{t})_{t \in [0,1]}$ defines a pinching deformation with support in $\mathcal{V}({\Lambda})$. Hence, we define $h_{t}$ the quasi conformal maps in $\hat{\C}$ given by the Measurable Riemann Mapping Theorem that integrates ${\sigma}_{t}$. Also, we can normalize and to assume that $h_{t}$ fixes $\infty$ and two points $p,q \in J(f)$. 

The composition $h_{t}\circ f \circ h_{t}^{-1}=f_{t}$ defines an entire transcendental map. We are interested in showing that pinching some Baker laminations we have that $f_{t} \rightrightarrows F$, where the sign $\rightrightarrows$ means uniform convergence, in this case $F$ is a meromorphic map with $\infty$ as essential singularity.

The following notions will be used in Lemma 2 and lemma 3.

Let us denote by $W_t$ the image of $W$ under $h_t$, if ${\gamma} \in W_t$ is a curve, denote by $l_{W_t}(\gamma)$ the hyperbolic length of $\gamma$ in $W$ and by $diam_{s}$ the spherical diameter of a set.

\begin{lema}Assume that $h_{t} \rightrightarrows H$, then for any $\gamma \in \mathcal{R}({\Lambda})$, $lim_{t \rightarrow 1}$ $(diam_{s}h_{t}(\bar{\gamma}))=0$
\end{lema}

\begin{proof}Notice that since the ellipse field is parallel to the geodesic $\gamma$, then $h_t(\gamma) \in W_t$, is a geodesic. Divide $\gamma$ in segments $s_i$ of the same length, then $l_{W_t}(s_i)< l_{W_r}(s_i)$ if $t<r$. Since $h_{t}$ converges,  $l_{W_t}(s_i)$ converges to zero as $t \rightarrow 1$, otherwise we could construct an ellipse field over the geodesic and we can  apply a pinching process again. Hence, the two points of $End(\gamma)$ must collide and we conclude that $lim_{t \rightarrow 1}$ $(diam_{s}h_{t}(\bar{\gamma}))=0$.

\end{proof}

\begin{lema} If $f$ satisfies conditions (a) and (b) above, the diameter of any sequence of elements in $\mathcal{V}({\Lambda})$ tends to $0$ if $\lambda_{\infty} \notin {\Lambda}$.
\end{lema}

\begin{proof} From conditions (a) and (b) at the beginning of the section, we have that $f$ is weakly hyperbolic, then the Ma\~{n}\'{e} shrinking lemma in [ST] can be applied to $f$ in $\C$. That implies that we have to avoid laminations involving $\infty$, by the discussion in subsection 3.2, such lamination contains $\lambda_{\infty}$. 
\end{proof}

\section{pinching along $\lambda_{\infty}$}

Consider $h_{t}$ a sequence of quasiconformal maps defined by pinching a Baker lamination as above.

\begin{lema}Let $f$ be a entire transcendental map with a univalent Baker domain $U$ of hyperbolic type I. Then, the pinching process along the lamination $\lambda_{\infty}$ does not converges. 
\end{lema}

\begin{proof} Let us assume that $f_{t} \rightrightarrows F$, by Lemma 2, $lim_{t \rightarrow 1}$ $(diam_{s}h_{t}(\overline{\lambda_{\infty}}))=0$ which implies that $h_{t}(\overline{\lambda_{\infty}}) \rightarrow \infty$. For $P \in \partial U$ the repulsive fixed point of $f$, let $P_{t}=h_{t}(P)$. Then $P_{t} \in \partial(h_{t}(U))$ is a  boundary point of $\lambda_{\infty}$ other than $\infty$ (along the process of pinching) and we have that $P_{t} \rightarrow {\infty}$, as $t \rightarrow 1$. Since $P_{t} \in J(f_{t})$, by the f-invariance of the pinching, we have that any point in the set $f_{t}^{-n}(P_{t})$ tends to $\infty$ under $h_{t}$, for any $n \geq 0$ and any $t<1$. However, the set $\cup_{n}f_{t}^{-n}(P_{t})$ is dense in $J(f_{t})$, then by uniform convergence we have that $J(F)=\infty$. That is a contradiction.

\end{proof}

As a direct consequence of this lemma, we have the following theorem.

\begin{theo}Let $f$ be a entire transcendental map with a periodic univalent Baker domain $U$ of hyperbolic type I. If $\Lambda$ is a Baker lamination on $U$ that contains $\lambda_{\infty}$, then, the pinching process along the lamination $\Lambda$ does not converges. 
\end{theo}

\section{pinching along general laminations}

Let us denote by $area(A)$ the plane area of a set $A$.

\begin{lema}$area(\mathcal{R}(\Lambda))=0$
\end{lema}

\begin{proof} $area({\Lambda})=0$ and for each $V \subset {\mathcal{U}}$ we have  $area({\mathcal {R}({\Lambda}})\bigcap V)=0$, then Lemma 3 implies $area(\mathcal{R}(\Lambda))=0$.
\end{proof}

\begin{defi}(McMullen [Mc1]) {\bf $J(f)$ is  thin at $\infty$} if for all $z\in \C$ there exist $R$ and $\epsilon$ such that $density(J(f), D_{r}(z))<1-\epsilon$, for all $r > R$.
\end{defi}

Where $D_{R}(z)$ is the disc centered at $z$ with radius $R$ and $density(E,G):=\frac{area(E \cap G)}{area G}$.

The main task of this section is to prove the following result:

\begin{theo}Let $f$ be an entire transcendental map satisfying the following properties:

(a) $f$ has a univalent Baker periodic domain $U$ of hyperbolic type I,

(b) the postcritical set $P(f)$ is a positive distance away from the Julia set.

(c) $J(f)$ is thin at $\infty$.

Then for  $\Lambda$ any Baker lamination in $U$ which does not contains $\lambda_{\infty}$, the pinching process along $\mathcal {R}(\Lambda)$ converges uniformly to an entire transcendental map $F$ which exhibits in its Fatou set a family of bounded simply connected wandering domains, disjoint from the postrcritial set.

\end{theo}

\begin{proof}
The detailed arguments are in [HT] and we will mention only the required results. By Ascoli's theorem, the sequence $\{h_{t}\}$ converges uniformly if and only if it is equicontinuous and converges pointwise.

The first step is to prove that the quasiconformal homeomorphisms $\{h_{t}\}$ are equicontinuous.
For that [HT] uses two lemmas: 

First Lemma: Equicontinuity criterion at a point, which is a criteria based on annuli neighborhoods at a point, and

Second Lemma: One good annulus around each Julia point in $\C$ (see Lemma 2.7 in [HP]).
This Lemma is proved using the construction in [HT] which in our case is a consequence of Lemma 1 and the Remark of this paper. Then the weak hyperbolicity condition on $f$, a consequence of hypothesis (b), is used to spread these annuli at every point in the Julia set. 

At $\infty$ we can always find a bounded annulus $A$ such that $\partial A \cap \mathcal{V}(\mathcal{R})=\emptyset$, because $\lambda_{\infty} \notin {\Lambda}$ and so the geodesics of the lamination either are contained in the annulus of are disjoint from it. Then there is $m >0$, such that $mod (h_{t}(A)) \geq m$ for all $t$. Satisfying the one good annulus criteria. Similarly it can be constructed a nested sequence of annuli $\{A_{n}\}$, with
$h_{t}(A_{n}) \geq m$.

Moreover, our Lemma 2 implies that when $h_{t_{n}} \rightrightarrows H$, the map $H$ sends each leaf of $\mathcal{R}$ to a point. By the discussion in Section 3.2, we have that this are the only fibers of $H$.

The second step is to prove that the limit map is unique.  

Let us prove first that $F$ is thin at infinity. Since the map $h_{1}=H$ is conformal and injective in ${\C}-{\mathcal{V}}(\overline{\Lambda})$, then the derivative of $H$ is non zero, therefore the one quarter Koebe theorem prevents for the components of      ${\mathcal{U}}-{\mathcal{V}}({\Lambda})$ to shrink excessively. From this situation and Lemma 5, we have that $J(f)$ and $J(F)$ have the same density for large discs $D_{r}(z)$. This implies that $J(F)$ is thin at infinity.

By Theorem B of Stallard [St], we conclude that the Lebesgue measure of the Julia set, $m(J(F))$, is $0$.

If we assume that there are two sequences $(t_{n})$ and $(s_{n})$ such tending to 1, such that
$h_{t_{n}} \rightrightarrows H_{1}$, $h_{s_{n}} \rightrightarrows H_{2}$, $f_{t_{n}} \rightrightarrows F_{1}$ and  $f_{s_{n}} \rightrightarrows F_{2}$, as explained in [HT] there is a homeomorphism $\varphi$ of $\hat{\C}$ conformal in $\hat{\C}-J(F_{1})$, such that $\varphi \circ F_{1} \circ \varphi^{-1}=F_{2}$. For weakly hyperbolic dynamics, it is proved in [H] (Proposition 6.3 and Theorem 6.1) that $\varphi$ is globally quasiconformal (his proof is general and it does not requires $f$ to be rational). Again, by [St], we have that $m(J(F_{1}))=0$, then $\varphi$ is a Moebius map. By normalization, we assume that our maps $h_{t}$ fixes three points, then $\varphi$ as well, consequently it is the identity, this implies that the limit is unique.

After pinching, by the contraction to a point of each leaf in the lamination, we obtain (from the Baker domain), a new domain in the Fatou set of $F$ that contains disjoint simply connected sets for which any point tends to $\infty$ under iteration. Then, we have a family of simply connected wandering domains disjoint from the postcritical set. This family can have one member or infinitely many members.

\end{proof}

The examples of Section 2.4, satisfies all the hypothesis of the theorem.

The topology of the family of wandering domains depends completely on the Baker lamination.

\begin{rema}
Pinching deformations can be done in the other types of Univalent Baker domains, but all known examples suggest that condition (b) is not satisfied. In such case we don't know that $m(J(F))=0$ and we can not assure the uniqueness of  the convergence.

\end{rema}

{\bf Question:} According to [R], Baker domains are in the low escaping sets, which implies that our wandering domains are also in the low escaping sets, however the doubly connected wandering domains in [KS] are disjoint from the postcritical set and in the fast escaping set. Bergweiler [B3], has an entire transcendental map with a fast escaping simply connected wandering domain, but with part of the postcritical set on it. Hence a natural question appear: Does there exist a simple connected wandering domain 
belonging to the fast escaping set but disjoint to the postcritical set?

\subsection{Shapes of the limit domains}\

In this section we will see that if the limit map $F$ exist as in theorem 2 above, then two things can occur. One is that the Baker domain ends up in a wandering domain and second that the Baker domain ends up into a wandering domain and a Baker domain. We can classify the laminations according to such phenomena in the following way:

 The closure of the geodesic $\lambda_{\infty}$, separates $\partial U$ into the left and right sides.

A) All the leaves of the lamination $\Lambda$ are in one side or in the other one.

B) Some leaf of the lamination $\lambda$ intersects $\lambda_{\infty}$.

In case (A), the connected component of $U-{\Lambda}$ to which $\lambda_{\infty}$ belongs,
develops in the limit of the deformation a univalent Baker domain. The other components will be wandering domains.

This observation proves the following:


\begin{theo} If $\Lambda$ is as in (A), then $U$ becomes in the limit of the deformation a Baker domain with a family of wandering domains attached to its boundary. If $\Lambda$
is as in (B), then, only a family of wandering domains appears.
\end{theo} 

\begin{figure}[h]
\begin{center}
\includegraphics[scale=0.5]{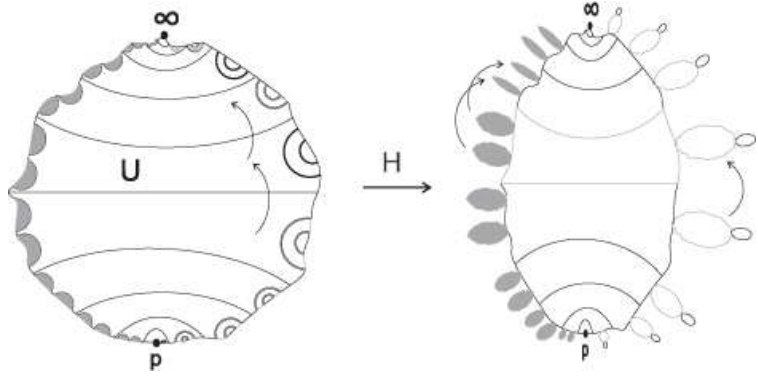}
\caption{A Baker domain $U$ with a Baker lamination (case A) and the result of the pinching: a new Baker domain and four different wandering domains.} 
\end{center}
\end{figure}

\section{Teichm\"uller space and Hurwitz class}

For the set of quasiconformal homeomorphisms of $\C$, define an equivalence relation $\sim$ by identifying $\varphi$ and $\psi$ if there exist a conformal map $c$, such that $c \circ\varphi=\psi$.

Denote by $QC(f)$ the set of quasiconformal homeomorphisms of $\C$ that commute with $f$ and by $QC_{0}(f)$ the subgroup of those quasiconformal homeomorphisms that are isotopic to the identity by an isotopy that fixes $Sing(f^{-1})$ pointwise.

Then define $Top(f)=\{\varphi:{\C}\rightarrow{\C}$ homeomorphism: there exist ${\psi}$ homeomorphism with ${\varphi}\circ f \circ {\psi}^{-1}=g$ an entire transcendental map $\}/{\sim}$. 

The group $QC(f)$ acts on $Top(f)$ by $(q,{\varphi})={\varphi}\circ q^{-1}.$

\begin{defi}
The {\bf Hurwitz marked class} ${\mathcal{M}}_{f}$ of an entire transcendental map $f$, is the space $Top(f)$ modulo the action of $QC_{0}(f)$, i.e. ${\mathcal{M}}_{f}=Top(f)/QC_{0}(f).$
\end{defi}

If $Sing(f^{-1})$ is discrete, $\varphi \in Top(f)$ can be choosen quasiconformal and in this case ${\psi}$ is also quasiconformal.

Eremenko and Lyubich proves in [EL] that if ${\phi}_{0} \circ f =g_{0} \circ {\psi}_{0}$ and ${\phi}_{1} \circ f =g_{1} \circ {\psi}_{1}$ and there is an isotopy ${\phi}_{t}$ connecting ${\phi}_{0}$ with ${\phi}_{1}$
fixing the singular values, then $g_{0}=g_{1}$.

Now, consider the surface $S={\C}-Sing(f^{-1})$  and let us denote by $T(S)$ the usual 
Teichm\"uller space of $S$. Given a quasiconformal map ${\varphi}\in{\mathcal{M}}_{f}$, the restriction of ${\varphi}$ to the surface $S$ gives an element in $T(S)$. Inversely, assuming that the measure of $Sing(f^{-1})$ is zero, given a quasiconformal map $\varphi$ in $T(S)$, it can be extend  to a quasiconformal map of $\C$. As shown in [EL] there is a homeomorphism $\psi$ and an entire transcendental map such that ${\varphi}\circ f=g \circ \psi$. This implies that there is a well defined bijective map from ${\mathcal{M}}_{f}$ to $T(S)$ which induces a topology in ${\mathcal{M}}_{f}$. 

For a definition of the Teichm\"uller space of $f$, $Teich(f)$ see for instance [McSu], [MNTU], [FH1]. According to the situation we describe above,  $Teich(f)$ is included in $T(S)$ as a subspace, hence $Teich(f)$ includes naturally as a subspace of ${\mathcal{M}}_{f}$.

Entire Transcendental maps satisfy properties (a) and (b) of Section 4, satisfies that  $m(J(f))=0$. According to [FH], we have $Teich(f)=Teich(U/f)\times_{i}Teich(D_{i}/f)$, being $D_{i}$ the full orbit of an invariant Fatou domains other than $U$ (there are no deformations in the Julia set).

For $f$ with a univalent periodic Baker domain of hyperbolic type I, we have that $Teich(f)$ is infinite dimensional complex Banach space . Infinite dimensional Teichm\"uller spaces have a Teichm\"uller metric which is non-differentiable and has arbitrarily short simple geodesics [EaLi], but are embedded into certain spaces of quadratic differentials by a theorem of Bers [GH]. Many other properties can be found in [NV], [GH] and [GL].

We are interested in the boundary of $Teich(f)$ as subspace of $M_{f}$.  For instance a consequence of Theorem 1 is that $Teich(f)$ is not compactly contained in $M_{f}$. By Theorem 2, the maps which are limits of pinching along Baker laminations other than $\lambda_{\infty}$ are in the boundary of $Teich(f)$ and are in $M_{f}$. These limits are similar to the cusp groups in the Bers boundary, whose density was proved by McMullen [Mc3] and are also reminiscent of the conjecture in [Mc2], that parabolic periodic domains are dense in the boundary.
Although no much is known about the boundary of Teichm\"uller spaces of infinite dimensions, it would be interesting to explore it and a problem is if our pinched limits are dense in the boundary of the slice $Teich(U/f)$. 

In the next proposition we are interested in the development of $U$ to a periodic domain as $(f_{t})$ tends to the boundary of $Teich(f)$, when $f_{t}$ deforms in the slice $Teich(U/f)$, this deformation is not necessarily by pinching a lamination.

\begin{prop}Assume that there is a sequence $(f_{t}) \in Teich(U/f)$ such that $(f_{t}) \rightrightarrows F$ then, if there are new Fatou periodic domains of $F$, they are univalent Baker periodic domains. 
\end{prop}

\begin{proof} Since $F \in M_{f}$ then $F$ is entire transcendental.  If $F$ has new Fatou periodic domains, then by the uniform convergence of $(f_{t})$ to $F$, it has to be a periodic univalent domain without singularities, hence a univalent periodic Baker domain.
The singular set of $F$  is contained in the persistent components of the Fatou domain.
\end{proof}

{\sl Acknowledgments}: To Peter Makienko who kindly explain us his early work on pinching, to Gwyneth Stallard for bringing our attention to her paper. To Memox for the figures. To PAPIIT grant IN105306.

\section{References}

[Ba] Baker I. N.: The domains of normality of an entire function. Ann. Acad. Sci. Fenn. A {\bf 1} (1975) pp. 277-283.

[B1] Bergweiler W.: Iterations of meromorphic functions. Bull. of the Amer. Math. Soc. (New Series). Vol 29, Number 2, October 1993.

[B2] Bergweiler W.: Invariant domains and singularities. Math. Proc. Camb. Phil. Soc. {\bf 117}
525-32.

[B3] Bergweiler W: An entire function with simply and multiply connected wandering domain. Pure and Applied mathematical Quarterly, Vol 7, number 1 (Special issue: In honor of Frederick Gehring, Part 1 of 2)(2011) pp:107-120.

[BE] Bergweiler W. and Eremenko A.: Direct singularities and completely invariant domains of entire functions. Illinois Jour. of Maths. 52 (2008) pp. 243-259.

[BF] Baransky K. and Fagella N.: Univalent Baker domains. Nonlinearity 14 (2001). Institute of Physics Publishing. pp 411-429.

[BH] Bullet S. and Ha\"{i}ssinsky P.: Pinching holomorphic correspondences. AMS Journal of Conformal Geometry and Dynamics {\bf 11} (2007) pp. 65-89.

[EaLi] Earle C. and Li Z.: Isometrically Embedded Polydiscs in Infinite Dimensional Teichm\"uller Spaces. The Journal of Geometric Analysis. Vol 9, Number 1, (1999), pp.51-71.

[EL] Eremenko A. and Lyubich M.: Dynamical properties of some classes of entire functions. Annales de l'institut Fourier, tome 42, no.4 (1992), pp 989-1020.

[EL2] Eremenko A. and Lyubich M.: Iterates of entire functions, Dokl. Akad. Nauk SSSR 279 (1984), 25-27; English transl. in Soviet Math. Dokl. 30 (1984).

[F] Fatou P. : Sur l'it{\'e}rations des fonctions transcendentes enti{\'e}res, Acta Math. {\bf 47} (1926) pp. 337-370.

[FH] Fagella N. and Hendriksen C.: The Teichm\"uller space of an entire function. Complex Dynamics: Families and Friends, AK Peters, Ltd. ISBN: 978-1-56881-450-6 

[FH1] Fagella N. and Hendriksen C.: Deformations of entire functions with baker domains, Discrete and Continuous Dynamical Systems 15 (2006), 379-394.

[GL] Gardiner F.P. and Lakic N.: Quasiconformal Teichm\"uller theory. Mathematical surveys and Monographs, vol. 76, Amer. Math. Soc., Providence, RI, 2000,xix+372pp.

[GH] Gardiner F.P. and  Harvey W.: Universal Teichm\"uller space. Handbook of Complex Analysis 1. Reiner Kuhnau ed. North-Holand, 2002.

[H] Ha\"{i}ssinsky P.: Rigidity and expansion for rational dynamics. J. London Math. Soc. (2) 63 (2001) pp. 128-140.

[He] Herman M.R.: Are there critical points on the boundaries of singular domains. Comm. Math. Phys 99 (1985), 593-612.

[HT] Haissisky P. and Tan Lei: Convergence of pinching deformations and matings of geometrically finite polynomials. Fund. Math. {\bf 181} (2004), p.143-188.

[HaTa] Harada T. and Taniguchi M.: On Teichm\"uller space of complex dynamics by entire functions. Bull. Hong Kong Math. Soc. 1 (1997), 257-266.

[KS] M. Kisaka and M. Shishikura: On multipliy connected wandering domains of entire functions, in "Transcendental Dynamics and Complex Analysis" (Ed. by P. J. Rippon and G. M. Stallard), London Mathematical Scociety Lecture Nore Series 348, Cambridge Univ. Press, 348 (2008), pp. 217-250.

[M] Makienko P.: Unbounded components in parameter space of rational maps, in: Conformal Geometry and Dynamics. Vol 4, 2000, p.1-21.

[Ma] Maskit B.: Parabolic elements in Kleinian groups. Ann. of Maths. 117 (1983)p.659-668.

[Mc1] McMullen C.: Area and Hausdorff dimension of Julia sets of Entire Functions. Trans. Amer. Math. Soc. Vol 300, Number 1, March 1987. pp. 329-342.

[Mc2] McMullen C.: Rational maps and Teichm\"uller spaces. Analogies and open problems. Lecture Notes in Mathematics. Vol. 1574. (1994). pp 430-433. Springer-Verlag. 

[Mc3] McMullen C.: Cusps are dense. Annals of Mathematics. 133 (1991) pp. 217-247.

[McSu] McMullen C and Sullivan D.: Quasiconformal homeomorphisms and dynamics III: The Teichm\"uller space of a holomorphic dynamica systems, Adv.Math.135 (1998),p.351-395

[MNTU] Morosawa S., Nishiura Y., Taniguchi M., Ueda T.: Holomorphic dynamics. Cambridge Studies in Advanced Mathematics 66. Cambridge University Press, 2000.

[NV] Nag S. and Verjovsky A.: $Diff(S^{1})$ and the Teichm\"uller spaces. Commun. Math. Phys. 130 (1990) pp. 123-138.

[R] Rippon P.: Baker domains of meromorphic functions. Erg. Theo. Dynam. Systems, 26 (2006), pp.1225-1233.

[RS] Rippon P. and Stallard G.: Families of Baker domains II. Conformal Geometry and Dynamics. Vol 3 (1999), pp. 67-78.

[Si] Sienra G.: Surgery and hyperbolic univalent Baker domains. Nonlinearity 19 (2006), pp.959-967.

[ST] Shishikura M. and L. Tan.: An Alternative proof of Ma\~{n}\'{e}'s theorem on non-expanding Julia sets. The Mandelbrot set, theme and variations, London Mathematical Society Lecture Notes Series 274 (ed. L. Tan, Cambridge University Press, 2000) pp.265-279.

[St] Stallard G.M.: Entire functions with Julia sets of zero measure.  Math. Proc. Cambridge Philos. Soc.,  vol. 108 issue 3 (1990) pp. 551--557.

[T] Tan Lei: On pinching deformations of rational maps. Ann. Sci.Ec.Nor.Sup. $4^{e}$ series,t.35, 2002, p.353-370.

\section{Keywords:} Entire transcendental maps, univalent Baker domains, laminations, pinching deformation, wandering domains, Teichm\"uller space, Hurwitz space.

{\bf Patricia Dom{\'i}nguez}, 

F.C. F{\'i}sico Matem{\'a}ticas.

Benem{\'e}rita Universidad Aut{\'o}noma de Puebla, M{\'e}xico.

Puebla, M{\'e}xico

pdsoto@fcfm.buap.mx

{\bf Guillermo Sienra}, 

Depto. de Matem{\'a}ticas. 

Fac. de Ciencias. UNAM. 

M{\'e}xico D.F., M{\'e}xico.

guillermo.sienra@gmail.com

\end{document}